\input amstex
\documentstyle{amsppt}
\pageheight{204mm}
\pagewidth{133mm}
\magnification\magstep1


\def\nologo{\let\logo@\empty}

\def\Ad{\operatorname{Ad}}

\def\Aut{\operatorname{Aut}}
\def\Aut{\operatorname{Aut}}

\def\Cont{\operatorname{Cont}}

\def\End{\operatorname{End}}
\def\Ext{\operatorname{Ext}}

\def\gr{\operatorname{gr}}

\def\Hom{\operatorname{Hom}}

\def\Im{\operatorname{Im}}

\def\Ker{\operatorname{Ker}}

\def\Mor{\operatorname{Mor}}

\def\rank{\operatorname{rank}}

\def\SL{\operatorname{SL}}

\def\Sym{\operatorname{Sym}}

\def\ts{\tsize\sum}

\def\bC{\bold C}

\def\bQ{\bold Q}

\def\bR{\bold R}
\def\bS{\bold S}

\def\bZ{\bold Z}

\def\cB{{\Cal B}}

\def\cF{{\Cal F}}

\def\cL{{\Cal L}}

\def\cO{{\Cal O}}

\def\cS{{\Cal S}}

\def\fa{{\frak a}}

\def\fg{{\frak g}}

\def\g{\gamma}
\def\G{\Gamma}

\def\l{\ell}
\def\sig{\sigma}
\def\Sig{\Sigma}

\def\vf{\varphi}

\def\.{$.\;$}

\def\fs{{\text{\rm fs}}}
\def\gp{{\text{\rm gp}}}
\def\loga{{\text{\rm log}}}

\def\resp.{\text{\rm resp}.\;}

\def\tra{\overset{\sim}\to{\to}}

\def\O^logten{\cO\^log\otimes}

\let\bs=\backslash

\let\lan=\langle
\let\lan=\langle

\let\lra=\longrightarrow

\let\hra=\hookrightarrow

\let\ox=\otimes

\let\ran=\rangle
\let\Ra=\Rightarrow

\let\sub=\subset
\let\sups=\supset

\let\x=\times

\def\Fb{\overline{F}}

\def\xb{\overline{x}}

\def\Dc{\check{D}}

\def\cExt{{\Cal E}xt}

\topmatter

\title
Log intermediate Jacobians
\endtitle

\author
Kazuya Kato, Chikara Nakayama, Sampei Usui
\endauthor

\address
\newline
{\rm Kazuya KATO}
\newline
Graduate School of Science
\newline
Kyoto University
\newline
Sakyo-ku, Kyoto, 606-8502, Japan
\newline
{\tt kzkt\@math.kyoto-u.ac.jp}
\endaddress

\address
\newline
{\rm Chikara NAKAYAMA}
\newline
Graduate School of Science and Engineering
\newline
Tokyo Institute of Technology 
\newline
Meguro-ku, Tokyo, 152-8551, Japan
\newline
{\tt cnakayam\@math.titech.ac.jp}
\endaddress

\address
\newline
{\rm Sampei USUI}
\newline
Graduate School of Science
\newline
Osaka University
\newline
Toyonaka, Osaka, 560-0043, Japan
\newline
{\tt usui\@math.sci.osaka-u.ac.jp}
\endaddress

\endtopmatter

\document

\head
Contents
\endhead

\S0. Introduction

\S1. Situation

\S2. Associated log Hodge structure

\S3. Notation

\S4. Main theorem

\S5. Example of a relatively complete fan

\S6. LMHS and their moduli

\S7. Log intermediate Jacobian

\S8. Relationship with other works

\S9. Example of slit

\S10. Comments on relative completeness

\S11. Case of extensions with negative weights

\S12. Remarks on N\'eron models 

\vskip30pt

\head
\S0. Introduction
\endhead

Let $S=\Delta=\{q\in \bC\;|\;|q|<1\}$, and let
$$
f:X\to S
$$ 
be a projective morphism which is smooth over
$S^*=\Delta^*=\Delta\setminus\{0\}$ and which is
semi-stable at $q=0$.

Put $X^*=f^{-1}(S^*)$.
Fix a polarization of $X^*/S^*$ and fix $r\ge1$.
Then we have the variation of polarized Hodge structure 
$H^{2r-1}(X^*/S^*)$ over $S^*$ (cf\.[Gri70]):
$$
\align
&H^{2r-1}_\bZ(X^*/S^*)
=R^{2r-1}f_*\bZ/(\text{torsion}),\\
&H_\cO^{2r-1}(X^*/S^*)
=R^{2r-1}f_*(\Omega^\bullet_{X^*/S^*})\simeq
\cO_{S^*}\ox_\bZ H_\bZ^{2r-1}(X^*/S^*),\\
&F^pH_\cO^{2r-1}(X^*/S^*)
=R^{2r-1}f_*(\Omega^{\ge p}_{X^*/S^*}).
\endalign
$$
Let us twist it into $H^{2r-1}(X^*/S^*)(r)$ whose weight is 
$-1$.
Fix a point $s\in S^*$.

Put
$$
H'=H^{2r-1}_{\bZ}(X^*/S^*)(r)_s, \qquad H=H' \oplus \bZ.
$$
Let $W$ be an
increasing filtration on $H_{\bQ}$ such that
$$
\gr^W_0(H_{\bQ})=\bQ,\quad
\gr^W_{-1}(H_{\bQ})=H'_{\bQ},\quad
\gr^W_k(H_{\bQ})=0\quad\text{for $k\ne0,-1$}.
$$
Let
$$
\align
&\lan\;,\;\ran_0:\bZ\x\bZ\to\bQ;\; (a,b)\mapsto ab,\\
&\lan\;,\;\ran'=\lan\;,\;\ran_{-1}:H'\x H'\to\bQ\quad
\text{the pairing defined by the polarization of 
$X^*/S^*$.}
\endalign
$$
Let $(h_k^{p,q})$ be the evident Hodge numbers.

Let $D'=D'(H',(h^{p,q}_{-1}),\lan\;,\;\ran')$ be the 
classifying space of polarized Hodge structures of 
weight $-1$, and let $D=D(H,W,(h^{p,q}_k),
(\lan\;,\;\ran_k))$ be the 
classifying space of gradedly polarized mixed Hodge 
structures.
Let $\gamma':H'\to H'$ be the local monodromy and let
$\Gamma'$ be the group generated by $\gamma'$.
Let $\Gamma$ be the subgroup of
$G_\bZ=\Aut(H,W,(\lan\;,\;\ran_k))$ consisting of all 
elements whose restrictions to $H'$ are contained in 
$\Gamma'$.

Then, the intermediate Jacobian $J^r_{X^*/S^*}$
introduced by Griffiths (cf\.[GS69]) is 
$$
J^r_{X^*/S^*}=(H_\cO^{2r-1}(X^*/S^*)
/(F^rH_\cO^{2r-1}(X^*/S^*)+H_\bZ^{2r-1}(X^*/S^*)))(r).
\tag0.1
$$
There is an isomorphism of functors (cf\.[Car80])
$$
\Mor(?,J^r_{X^*/S^*})\simeq
\cExt^1_{MHS}(\bZ,H^{2r-1}(X^*/S^*)(r))
\tag0.2
$$
from the category of complex spaces over $S^*$ to the
category of sets.
Here $\cExt^1$ is taken in the category (MHS) of
mixed Hodge structures with polarized graded
quotients (= Variation of gradedly polarized mixed
Hodge structure without assuming Griffiths
transversality).
This is also expressed as the fiber product (cf\.[Usu84])
$$
\CD
J^r_{X^*/S^*}@>>>\Gamma\bs D\\
@VVV@VV\gr^W_{-1}V\\
S^*@>\vf>>\;\;\Gamma'\bs D',
\endCD
\tag0.3
$$
where $\vf$ is the period map associated to the VPHS
$H^{2r-1}(X^*/S^*)(r)$.

The purpose of this article is to extend $J^r_{X^*/S^*}$
in expression (0.3) to a \lq\lq log intermediate
Jacobian $J^r_{X/S,\Sig}$" over $S$ by using 
log geometry (\S7).
Our \lq\lq$J^r_{X/S,\Sig}$" is related to 
works by Zucker [Zuc76], Clemens [Cle83], and Saito
[Sai96] (\S8).
 
\bigskip

  \S\S1--10 in this article are based on 
our note on which the third author gave a 
talk 
at the JAMI conference at Johns Hopkins university in March, 2005. 

\bigskip

\head
\S1. Situation
\endhead

Let $S=\Delta=\{q\in \bC\;|\;|q|<1\}$, and let
$$
f:X\to S
$$ 
be a projective morphism which is smooth over
$S^*=\Delta\setminus\{0\}$ and which is semi-stable
at $q=0$.
Endow $X$ and $S$ with the standard log structures
([Kak89]):
$$
\aligned
&M_S:=\{g\in\cO_S\,|\,\text{$g$ is invertible
outside the point $0\in S$}\}\hra\cO_S,\\
&M_X:=\{g\in\cO_X\,|\,\text{$g$ is invertible
outside the fiber $f^{-1}(0)$}\}\hra\cO_X.
\endaligned
$$
Associated with these, [KN99] introduced ringed spaces 
$(S^\loga,\cO_S^\loga)$, $(X^\loga,\cO_X^\loga)$, 
a morphism $f^\loga:X^\loga\to S^\loga$, and the 
following commutative diagram:
$$
\CD
X@<\tau_X<<X^\loga\\
@VfVV@VV{f^\loga}V\\
S@<\tau_S<<\;S^\loga.
\endCD
$$

An intrinsic definition of them is as follows.
As a set, $X^\loga$ is defined to be the set of all
pairs $(x,h)$ consisting of a point $x\in X$ and an
{\it argument function} $h$ which is a
homomorphism $M_{X,x}\to\bS^1$ whose restriction
to $\cO_{X,x}^\x$ is $u\mapsto u(x)/|u(x)|$.
Here $\bS^1:=\{z\in\bC\;|\;|z|=1\}$.
For definition of the topology of $X^\loga$, we work
locally on $X$.
Take a chart $\cS\to M_X$, then we have an injective
map
$$
X^\loga\hra X\x\Hom(\cS^\gp,\bS^1),\;\;
(x,h)\mapsto(x,h_\cS),
$$
where $h_\cS$ denotes the composite map
$\cS^\gp\to M_{X,x}^\gp\overset h\to\to\bS^1$.
We endow $X^\loga$ with the topology as a subset
of $X\x\Hom(\cS^\gp,\bS^1)$.
This topology is independent of the choice of a chart,
and hence is globally well-defined.
The canonical map
$$
\tau_X:X^\loga\to X,\;\;(x,h)\mapsto x,
$$
is surjective, continuous and proper.
For $x\in X$, the inverse image $\tau^{-1}(x)$ is
homeomorphic to $(\bS^1)^r$, where
$r:=\rank_\bZ(M_X^\gp/\cO_X^\x)_x$.
The continuous map $\tau_S:S^\loga\to S$ is defined
similarly, and the continuous map $f^\loga:X^\loga\to
S^\loga$ is defined to be $(x,h)\mapsto(f(x),
M_{S,f(x)}@>f^*>>M_{X,x}@>h>>\bS^1)$.

The sheaf of rings $\cO_X^\loga$ on $X^\loga$ is
defined as follows.
We will define first the {\it sheaf of logarithms}
$\cL$ of $M_X^\gp$ on $X^\loga$, and will then
define $\cO_X^\loga$ as a sheaf of
$\tau^{-1}(\cO_X)$-algebras generated by $\cL$, where
$\tau=\tau_X$.
Let $\cL$ be the fiber product of
$$
\CD
@.\tau^{-1}(M_X^\gp)\\
@.@VVV\\
\Cont(\quad,i\bR)@>{\exp}>>\;\Cont(\quad,\bS^1),
\endCD
$$
where $\Cont(\quad,T)$, for a topological space
$T$, denotes the sheaf on $X^\loga$ of continuous
maps to $T$, and $\tau^{-1}(M_X^\gp)\to
\Cont(\quad,\bS^1)$ comes from the definition of
$X^\loga$.
We define
$$
\cO_X^\loga:=
(\tau^{-1}(\cO_X)\ox_\bZ\Sym_\bZ(\cL))/\fa,
$$
where $\Sym_\bZ(\cL)$ denotes the symmetric
algebra of $\cL$ over $\bZ$, and $\fa$ is the ideal
of $\tau^{-1}(\cO_X)\ox_\bZ\Sym_\bZ(\cL)$
generated by the image of
$$
\tau^{-1}(\cO_X)\to
\tau^{-1}(\cO_X)\ox_\bZ\Sym_\bZ (\cL),\quad
f\mapsto f\ox1-1\ox\iota(f).
$$
Here the map $\iota:\tau^{-1}(\cO_X)\to\cL$ is the
one induced by
$$
\align
&\tau^{-1}(\cO_X)\to\Cont(\quad,i\bR),\quad
f\mapsto\tfrac12(f-\bar f),\quad\text{and}\\
&\tau^{-1}(\cO_X)\overset\exp\to\lra
\tau^{-1}(\cO_X^\x)\sub\tau^{-1}(M_X^\gp).
\endalign
$$
In the above, $\,\bar{}\;$ means the complex
conjugation.
We denote the projection $\cL\to
\tau^{-1}(M_X^\gp)$ by $\exp$, and the inverse
$\tau^{-1}(M_X^\gp)\to\cL/(2\pi i\bZ)$ by $\log$.
Then we have a commutative diagram with exact
rows:
$$
\CD
0@>>>\bZ@>2\pi i>>\tau^{-1}(\cO_X)
@>{\exp}>>\tau^{-1}(\cO_X^\x)@>>>1\;\\
@.@|@V{\iota}V{\cap}V@VV{\cap}V@.\\
0@>>>\bZ@>2\pi i>>\cL@>{\exp}>>
\tau^{-1}(M_X^\gp)@>>>1.
\endCD
$$

The morphism $f:X\to S$ induces a morphism
$$
f^\loga:(X^\loga,\cO_X^\loga)
\to(S^\loga,\cO_S^\loga)
$$
of ringed spaces over $\bC$ in the evident way.

[Usu01] showed that $f^\loga:X^\loga\to S^\loga$ is
topologically locally trivial over the base.

For $y\in X^\loga$, the stalk $\cO_{X,y}^\loga$ is
described as follows.
Put $x=\tau(y)\in X$ and
$r=\rank_\bZ(M_X^\gp/\cO_X^\x)_x$.
Let $(\l_j)_{1\le j\le r}$ be a family of elements
of $\cL_y$ whose images in
$(M_X^\gp/\cO_X^\x)_x$ form a system of free
generators.
Then we have an isomorphism of
$\cO_{X,x}$-algebras
$$
\tau^{-1}(\cO_{X,x})[T_1,\dots,T_r]
\tra\cO_{X,y}^\loga,\;\;T_j\mapsto\l_j.
$$
Note that this is a polynomial algebra of $r$-variables
over $\tau^{-1}(\cO_{X,x})$, which is {\it not} a local
ring if $r\ge1$.
\medskip

\head
\S2. Associated log Hodge structure
\endhead

For each $m\ge0$, we have the polarizable log Hodge
structure $H^m(X/S)$ of weight $m$:
$$
H^m_\bZ(X/S)=R^mf^\loga_*\bZ/(\text{torsion}),\quad
H_\cO^m(X/S)=R^mf_*(\Omega^\bullet_{X/S}(\log)),
$$
where $\Omega_{X/S}^\bullet(\log)$ is the de Rham
complex with log poles, with the Hodge filtration 
$$
F^pH_\cO^m(X/S)=R^mf_*(\Omega^{\ge p}_{X/S}(\log))
$$
and with the isomorphism
$$
\iota: \cO_S^{\log}\ox_\bZ H_\bZ^m(X/S)\simeq\cO_S^\loga
\ox_\cO H_\cO^m(X/S)
$$
extending the isomorphism
$\cO_{S^*}\ox_\bZ H_\bZ^m(X/S)|_{S^*}\simeq
H_\cO^m(X/S)|_{S^*}$, which satisfy, for any point
$t\in\tau_S^{-1}(0)$ and any ring homomorphism
$a:\cO_{S,t}^\loga\to\bC$ extending $\cO_{S,0}\to\bC$;
$g\mapsto g(0)$ with $\Im((2\pi i)^{-1}a(\log q))\gg0$,
$(H^m_\bZ(X/S)_t,F(a))$ is a plolarizable Hodge
structure ([Kaf98], [Mat98], [KU99], [KMN02], [KU09]).
  Here $F(a)$ is the filtration induced by the 
isomorphism $\iota$ and the ring homomorphism $a$. 

This is a reformulation, in terms of log Hodge theory,
of the classical theory of Schmid [Sch73], and Steenbrink
[Ste76] which says that a limit Hodge structure appears
at $q=0$.
\medskip

\head
\S3. Notation
\endhead

Fix a polarization of $X/S$.
Fix $r\ge1$, and consider the polarized log Hodge
structure $H^{2r-1}(X/S)(r)$ of weight $-1$ over $S$.
Fix a point $s\in S^*=\Delta^*$.
Put
$$
H'=H^{2r-1}_{\bZ}(X^*/S^*)(r)_s, \qquad H=H' \oplus \bZ.
$$
Let $W$ be an
increasing filtration on $H_{\bQ}$ such that
$$
\gr^W_0(H_{\bQ})=\bQ,\quad
\gr^W_{-1}(H_{\bQ})=H'_{\bQ},\quad
\gr^W_k(H_{\bQ})=0\quad\text{for $k\ne0,-1$}.
$$
Let
$$
\align
&\lan\;,\;\ran_0:\bZ\x\bZ\to\bQ;\; (a,b)\mapsto ab,\\
&\lan\;,\;\ran'=\lan\;,\;\ran_{-1}:H'\x H'\to\bQ\quad
\text{the pairing defined by the polarization of 
$X/S$.}
\endalign
$$
Let $(h_k^{p,q})$ be the evident Hodge numbers.
\medskip

\head
\S4. Main theorem
\endhead

Let $\gamma':H'\to H'$ be the local monodromy, i.e., the action of the 
standard generator of $\pi_1(S^*)$, and let 
$\Gamma' \subset \Aut(H', \lan\;,\;\ran')$ be the group generated by $\gamma'$.
Let $\Gamma$ be the subgroup of
$G_\bZ=\Aut(H,W,(\lan\;,\;\ran_k))$ consisting of
all elements whose restrictions to $H'$ are contained in
$\Gamma'$.

Let $N'=\log(\gamma'):H'_{\bQ}\to H'_{\bQ}$ (note that $\gamma'$ is
unipotent).
Let $\Sig'$ be the fan in
$\fg':=\End(H'_\bR,\lan\;,\;\ran') =\{h:H'_\bR\to
H'_\bR\;|\;\lan h(x), y\ran'+\lan x,h(y)\ran'=0\}$
defined by
$$
\Sig'=\{(\bR_{\ge0})N',\{0\}\}.
$$

\proclaim{Theorem}
There is a fan $\Sig$ in $\fg:=\End(H_\bR,W,
(\lan\;,\;\ran_k))$,
consisting of rational nilpotent cones, which satisfies the
following $(1)$ and $(2)$. 
\smallskip

$(1)$ For any $\sig\in\Sig$, $\sig$ is 
admissible for $W$ and the restriction of any element of $\sig$ to $H'_{\bR}$
is contained in $(\bR_{\ge0})N'$.
Furthermore, $\Sig$ is strongly compatible with $\G$.
\smallskip

$(2)$ $($Relative completeness.$)$ Let $\sig$ be any rational
nilpotent cone in $\fg$ which is
admissible for $W$ such that the restriction of
any element of $\sig$ to $H'_{\bR}$ is contained in $(\bR_{\ge0})N'$.
Then there exists a finite subdivison $\{\sig_j\}$ of
$\sig$ such that each $\sig_j$ is contained in some
element of $\Sig$.
\endproclaim

We explain some terminology in the above theorem.
A {\it nilpotent cone $\sig$ in $\fg$} is a cone over
$\bR_{\ge0}$ in
$\fg$ generated by a finite number of mutually
commutative nilpotent elements of $\fg$. We also
assume that $\sig$ is sharp, i.e., $\sig\cap(-\sig)=\{0\}$.
A nilpotent cone $\sig$ is {\it rational} if its generators
as an $\bR_{\ge0}$-cone can be taken in $\fg_\bQ$.
A nilpotent cone $\sig$ is {\it admissible for $W$} if,
for any element $N$ of $\sig$, there exists the
$W$-relative
$N$-filtration $M=M(N,W)$, i.e., the filtration
$M$ on $H_\bR$ satisfying
$$
\align
&NM_k\sub M_{k-2}\quad\text{for any }k,\qquad\text{and}\\
&N^l:\gr^M_{k+l}\gr^W_k\tra\gr^M_{k-l}\gr^W_k\quad
\text{for any }k, l\ge0.
\endalign
$$
Furthermore, this filtration $M$ depends only on the smallest face of 
$\sig$ containing $N$.


See [SZ85], [Kas85], [Kas86] for the details 
on admissibility.

See [KU09] for the definition of the strong compatibility.

\definition{Definition}
  We call a fan $\Sig$ in $\fg$ satisfying the conditions (1) and (2) 
of Theorem as a {\it relatively complete} fan. 
\enddefinition

\medskip

\head
\S5. Example of a relatively complete fan
\endhead

  We prove the main theorem by giving an explicit example of a 
relatively complete fan. 
  Let
$$
P:=\Im(N':H'_\bQ\to H'_\bQ),\quad
Q:=\Ker(N':H'_\bQ\to H'_\bQ)\cap P.
$$
  Take a finitely generated $\bZ$-submodule $L$ of $H'_{\bQ}$
containing $H' + N'(H')$.
Take a section  $s:(P\cap L)/(Q\cap L)\to P\cap L$ of $\bZ$-modules 
to the surjection 
$P\cap L \to (P\cap L)/(Q\cap L)$.
  Write as $s$ also for the induced $\bQ$-linear map
$P/Q\to P$ by abuse of notation. 
For an element $x$ of $P/Q$, let $a(x)$ be the order 
of the image of $x$ in $P/((P\cap L)+Q)$. 

Fix a $\bZ$-basis $(e_j)_{1\le j\le m}$ of $Q\cap L$, and
fix an element $e\in H$ inducing $1\in\bZ=\gr^W_0(H)$.

For $x\in P/Q$ and $n=(n_j)_j\in\bZ^m$, let
$\sig(x,n)$ be the $(\bR_{\ge0})$-cone in
$\fg$ generated by all elements $N\in\fg$ having the
following property:
\smallskip

The restriction of $N$ to $H'_{\bR}$ coincides with $N'$, and
$N(e)$ is an element of $H'_\bR$ of the form $s(x)+
(1/a(x))\sum_{1\le j\le m}c_je_j$ with $n_j\le c_j\le n_j+1$
for all $j$. 
\smallskip

Let $\Sig$ be the set of all faces of $\sig(x,n)$ for
all $x \in P/Q$ and $n\in\bZ^m$.

  It is easy to see that this $\Sig$ satisfies the desired conditions 
(1) and (2), 
and also readers can find a proof in \S11 in a more general setting. 
\medskip

\head
\S6. LMHS and their moduli
\endhead


In this section, we review the notion of LMHS and their
moduli in \cite{KNU08}, \cite{KNU09}, [KNU.p].
\medskip

A {\it log mixed Hodge structure with polarized graded
quotients} over an fs log analytic space $S$ is
$(H_\bZ,H_\cO,F,W,(\lan\;,\;\ran_k))$, consisting of a
locally constant sheaf of free 
$\bZ$-modules $H_\bZ$ on $S^\loga$, a locally free
$\cO_S$-module $H_\cO$ on $S$, with
$\cO_S^\loga\ox_\bZ H_\bZ\simeq
\cO_S^\loga\ox_{\cO_S}H_\cO$, 
a decreasing filtration $F$ on $H_\cO$ such that 
$\gr_F^p$ is a locally free $\cO_S$-module for any $p$,
a $\bQ$-rational increasing filtration $W$ of the sheaf $H_{\bR}$, 
and a family of polarizations $\lan\;,\;\ran_k$ on each $\gr_k^W$, 
which satisfies the following conditions (1)--(3)
(See [KU09] 2.6, \cite{KKN08a} 2.3, 2.5 for the details):

(1) {\it Admissibility} for 
any $s\in S$: the action of $\pi_1(s^{\log})$ (= the local monodromy 
action at $s$) is admissible 
with respect to $W$. 

(2) {\it Small Griffiths transversality} for any $\sig$. 
(This is a weaker version of
Griffiths transversality, and is imposed only at the place of degeneration.)

(3) {\it Graded positivity}: $(W_kH_\bZ/W_{k-1}H_\bZ,
F(\gr^W_k),\lan\;,\;\ran_k)$ is a polarized log Hodge structure of weight $k$ 
for each $k$.
\medskip

Let $D=D(H,W,(h_k^{p,q}),(\lan\;,\;\ran_k))$ be the
classifying space of log mixed Hodge structures with
polarized graded quotients, and let
$\Dc=\Dc(H,W,(h_k^{p,q}),(\lan\;,\;\ran_k))$ be its
\lq\lq compact dual".
They and their toroidal partial enlargements 
are defined roughly as follows.
  See [KU99], [KU09], [KNU08],
[KNU09], [KNU.p] for the details.

$$
\align
\cF&:=\{F\;\text{decreasing filtration on $H_\bC$}\;|\;
\dim_\bC\gr^p_F\gr^W_k=h_k^{p,w-p}\}\\
&\;\,\sups\Dc:=\{F\in\cF\;|\;\text{$\lan F^p(\gr^W_k),
F^q(\gr^W_k)\ran_k=0$ for $p+q>k$}\}\\
&\;\,\sups D:=\{F\in\Dc\;|\;\text{$i^{p-q}\lan x,
\xb\ran_k>0$ for $\;x\in H^{p,q}(\gr^W_k)-\{0\}$,
$p+q=k$}\},
\endalign
$$
where $H^{p,q}(\gr^W_k)=F^p(\gr^W_k)\cap
\Fb^q(\gr^W_k)$ for $p+q=k$.


  Let $\Sig$ be a fan satisfying (1) in Theorem. 

$D_\Sig=D_\Sig(H,W,(h_k^{p,q}),(\lan\;,\;\ran_k),\Sig)$
be the space of nilpotent orbits.

$\Dc'=\Dc'(H',(h_{-1}^{p,q}),\lan\;,\;\ran')\sups 
D'=D'(H',(h_{-1}^{p,q}),\lan\;,\;\ran')$

$D'_{\Sig'}=D'_{\Sig'}(H',(h_{-1}^{p,q}),\lan\;,\;\ran',\Sig')$
be the space of nilpotent orbits.

Then $D'=\gr^W_{-1}(D)$,
$\Dc'=\gr^W_{-1}(\Dc)$,
$D'_{\Sig'}=\gr^W_{-1}(D_\Sig)$.
\medskip

The polarized log Hodge structure
$H^{2r-1}(X/S)(r)$ defines a canonical morphism (period
map)
$$
S\to\Gamma'\bs D'_{\Sig'}.
$$

\head
\S7. Log intermediate Jacobian
\endhead

  Let $\Sig$ be a fan satisfying (1) in Theorem. 
  Then 

\proclaim{Proposition}
The space 
$\Gamma\bs D_{\Sig}$ is a log manifold 
and Hausdorff. 
\endproclaim

  This is a main result of \cite{KNU.p}, 
which is a mixed version of \cite{KU09}, and 
proved similarly to the pure case in loc.\ cit.  
  As mentioned in \cite{KNU08} 0.8, in the proof, 
the SL(2)-orbit theorem 0.5 in loc.\ cit.\ plays a key role, 
that is, to guarantee the continuity of the CKS-map 
(cf.\ \cite{KU09}), exactly 
as the SL(2)-orbit theorem of Cattani-Kaplan-Schmid [CKS86] 
did so in the pure case \cite{KU09}.

  Define
$J^r_{X/S,\Sig}$ to be the fiber product
$$
\CD
J^r_{X/S,\Sig}@>>>\Gamma\bs D_{\Sig}\\
@VVV@VV\gr^W_{-1}V\\
S@>\text{period map}>>\;\;\Gamma'\bs D'_{\Sig'}
\endCD
$$
in the category $\cB(\log)$ (\cite{KU09} 3.2.4).

\proclaim{Corollary}
  The space $J^r_{X/S, \Sig}$ is a log manifold and Hausdorff. 
\endproclaim

If $T$ also satisfies (1) in Theorem and if $T\sub\Sig$,
then $J^r_{X/S,T}$ is an open subset of
$J^r_{X/S,\Sig}$.

For $\Sig$ which also satisfies (2) in Theorem, i.e., 
for a relatively complete fan $\Sig$,  we call
$J^r_{X/S,\Sig}$ the {\it $r$-th log intermediate
Jacobian} associated to $X/S$ and $\Sig$.

\medskip

\head
\S8. Relationship with other works
\endhead

What we have just constructed is closely related to 
some works of Zucker [Zuc76] (cf\.also [EZ84]), Clemens
[Cle83], and Saito [Sai96]. 
What they considered are essentially
$J^r_{X/S,\Sig}$ for the following fans $\Sig=\Sig_j$
$(j=0,1,2)$.
We denote Zucker's (resp\.Clemens', Saito's)
fan by $\Sig_0$ (resp\.$\Sig_1$, $\Sig_2$).
We also denote our relatively complete fan in \S5 by $\Sig_3$.
$\Sig_j$ for $j=1,2$ are considered under the condition:
\medskip

\noindent
(8.1)\qquad
$(N')^2=0$ and $\gr^{W'}_0$ has Hodge type $(0,0)$.
\medskip

\noindent
Here $W'=W(N')[1]$, that is, the $N'$-filtration on $H'$.

Now we define $\Sig_j$ $(j=0,1,2)$:
$$
\align
&\Sig_0=\{(\bR_{\ge0})N\;|\;N\in\fg_\bQ,N|_{H'}=N',
N(e)\in N'(H')\}\cup\{\{0\}\}.\\
&\Sig_1=\{(\bR_{\ge0})N\;|\;N\in\fg_\bQ,N|_{H'}=N',
N(e)\in Q \cap H'\}\cup\{\{0\}\}.\\
&\Sig_2=\{\text{face of $\sig_n$}\;|\;n\in\bZ^m\},
\endalign
$$
where for a fixed $\bZ$-basis $e'_1,\dots,e'_m$ of
$N'(H')$, $\sig_n$ is the $(\bR_{\ge0})$-cone generated
by
$$
\{N\in\fg_\bQ\;|\;N|_{H'}=N',
\text{$N(e)=\ts_{1\le j\le m}c_je'_j$ with
$n_j\le c_j\le n_j+1$ for any $j$}\}.
$$

Under the condition (8.1), we see that $P=Q$ so 
that $s$ must be $0$,  and 
for any $L$, 
we have the following
relationship among the four fans:
$$
\matrix
\Sig_2&<&\Sig_3\\
\cup&&\cup\\
\Sig_0&\sub&\Sig_1\\
\endmatrix,\qquad
\text{$\Sig_1$ and $\Sig_2$ are not necessarily contained in each
other.}
$$
Here $A<B$ means $B$ is a subdivision of $A$.

Precisely speaking, as Saito pointed out in [Sai96,
(3.5) (iv)], Zucker's space is not Hausdorff.
This is because Zucker did not put the Griffiths
transversality.
On the other hand, 
our space $J^r_{X/S,\Sig_0}$ or $J^r_{X/S,\Sig_3}$ is always Hausdorff
thanks to  slits coming from the Griffiths transversality (see the next section).
In the case (8.1), the Griffiths transversality is
automatically satisfied and hence slits do not appear,
and Clemens and Saito considered exactly
$J^r_{X/S,\Sig_j}$ ($j=1,2$, respectively) which is
Hausdorff.
Clemens' is \lq\lq N\'eron model" which is not necessarily proper
over $S$, whereas Saito's is proper over $S$.

\proclaim{Proposition}
  Let $\Sig$ be a fan satisfying $(1)$ in Theorem. 
  Assume $(8.1)$. 
  Then,  $J^r_{X/S,\Sig} \to S$ is proper if and only if 
$\Sig$ is relatively complete. 
\endproclaim

  This is because slits do not appear under (8.1). 
\medskip

\head
\S9. Example of slit
\endhead

Let $Y=\bC/\bZ[i]$ be the elliptic curve with period $i$,
and let $g:E\to\Delta$ be the standard degeneration of
elliptic curves, i.e., $g^{-1}(q)=\bC^\x/q^\bZ$ for
$q\ne0$ and $g^{-1}(0)$ is a rational curve with one
node, and consider the family
$$
f:=g \circ {\text{pr}}_2:X=Y^2\x E\to S=\Delta.
$$
Then we see
$$
\text{(Zucker's space for $\Sig_0$)}=
Y^4\x(((\bC^\x)^4\x\bC^2\x\Delta)/\sim).
$$
Here $((t_j),(a_k),q)\sim((t'_j),(a'_k),q')$ if and only if
$$
\cases
\text{when $q\ne0$};\quad q'=q,\;
t'_j/t_j\in q^\bZ\;\text{for any }j,
\text {there exists }b\in\bZ[i]\;\\
\qquad\;\text{such that}\;
a'_2-a_2-b=0\;\text{and}\;
a'_1-a_1-b\cdot\log(q)/(2\pi i)\in\bZ[i];\\
\text{when $q=0$};\quad q'=q=0,\;
t'_j=t_j\;\text{for any }j,\;a'_2-a_2=0,\;a'_1-a_1\in\bZ[i].
\endcases
$$
This 
is not Hausdorff.
In fact, let $c\in\bC$ and $t\in(\bC^\x)^4$, and for
$n\gg0$, let $q_n=\exp(2\pi i(c+ni))$.
Then $(t,(c,1),q_n)\sim(t,(0,0),q_n)$, and $(t,(c,1),q_n)$
(resp\.$(t,(0,0),q_n)$) converges to $(t,(c,1),0)$
(resp\.$(t,(0,0),0)$) as $n\to\infty$.
But we see $(t,(c,1),0)\not\sim(t,(0,0),0)$.

On the other hand, our space for the fan $\Sig_0$ is
$$
J^2_{X/S,\Sig_0}=
Y^4\x\{\text{class of $(t,a,q)$}\;|\;q=0\Ra a_2=0\}.
$$
Here the condition \lq\lq$q=0\Ra a_2=0$" comes from
the Griffiths transversality, produces a slit and makes
the space $J^2_{X/S,\Sig_0}$ Hausdorff for the usual
topology and hence for the strong topology.

\medskip

\head
\S10. Comments on relative completeness
\endhead

We have an embedding
$$
\Mor(?,J^r_{X/S,\Sig})\sub
\cExt^1_{LMHS}(\bZ,H^{2r-1}(X/S)(r))
$$
of functors from the category $(\fs/S)$ of fs log
analytic spaces over $S$ to the category of
sets.
Here, $\cExt^1$ is taken in the category (LMHS) of log mixed
Hodge structures with polarized graded quotients (\S6).

\proclaim{Proposition} Let $\Sig$ be a relatively complete fan ({\rm \S4}). 
Then, for any fs log analytic space $S'$ over $S$ and any
$a\in\Ext^1_{S'}(\bZ,H^{2r-1}(X/S)(r))$, locally on $S'$,
there is a log modification $S''\to S'$ ({\rm {\cite{KU09}} 3.6})
and a subdivision $\Sig'$ of $\Sig$ satisfying $(1)$ in Theorem 
in \S$4$ 
such that the image of $a$ in
$\Ext^1_{S''}(\bZ,H^{2r-1}(X/S)(r))$ belongs to
$\Mor(S'',J^r_{X/S,\Sig'})$.
\endproclaim

  The proof of the above fact is similar to the pure case \cite{KU09} 4.3, 
where the extensions of period maps are explained. 

  In particular, we have the following. 
For any fs log analytic space $S'$ over $S$ which is log smooth over 
$\bC$ (\cite{KU09} 2.1.11), 
let 
$U$ be the open subspace of $S'$ where the log structure is trivial. 
  Let $a\in\Ext^1_{U}(\bZ,H^{2r-1}(X/S)(r))$ be an 
extension of graded polarized variation of MHS, 
regarded as a morphism $a: U \to J^r_{X_U/U}$ to the usual intermediate 
Jacobian. 
  Assume that $a$ is admissible with respect to 
$S'$. 
  Then, locally on $S'$,
there is a log modification $S''\to S'$ with $U \subset S''$ 
and a subdivision $\Sig'$ of $\Sig$ satisfying (1) in Theorem 
such that 
$U @>a>>  J^r_{X_U/U} \to J^r_{X/S,\Sig'}$ 
extends to a morphism 
$S'' \to J^r_{X/S,\Sig'}$. 

  More specifically,  assume $S'=S = \Delta$. 
  Then, $\Ext^1_{S'}(\bZ,H^{2r-1}(X/S)(r))$ is nothing but the space 
of admissible normal functions (\cite{Sai96}), 
and the above fact says that any 
admissible normal function extends to some log intermediate Jacobian 
because, in this case, there is no non-trivial log modification, 
that is, $S'' = S'$. 

  Since a cycle on $X$ gives an admissible normal function 
by a theory of Saito 
([Sai90], [Sai96]), we also have the Abel-Jacobi
map into the log intermediate Jacobian.

\bigskip

\head
\S11. Case of extensions with negative weights
\endhead

Here we give a construction of a
relatively complete fan for an extension 
$0\to(H'\text{ of weight $k$)}\to H \to\bZ\to0$
with $k<0$, which generalizes the case $k=-1$ in \S5. 
  Here the base $S$ is the disc $\Delta$.

\medskip

\noindent
{\bf {11.1.}} 
Let $W'$ be the filtration $W(N')[-k]$ on
$H'_\bQ$. Let
$$
P=\Im(N':H'_\bQ\to H'_\bQ)+W'_{-2},\quad
Q=\Ker(N':H'_\bQ\to H'_\bQ)\cap W'_{-2}.$$
Let $L$ be a finitely generated $\bZ$-submodule of $H_{\bQ}$ which 
contains $H' + N'(H')$. 
Fix a homomorphism $s:(P\cap L)/(Q\cap L)\to P\cap L$ of $\bZ$-modules  such that the composition $(P\cap L)/(Q\cap L) \overset s \to 
\to P\cap L \to (P\cap L)/(Q\cap L)$ is the identity map of $(P\cap L)/(Q\cap L)$, and denote the $\bQ$-linear map
$P/Q\to P$ induced by $s$ by the same letter $s$. 
For an element $x$ of $P/Q$, let $a(x)$ be the smallest integer $a$ such that $ax\in (P\cap L)/(Q\cap L)$. 

Fix $e\in H$ whose image in $H/H'=\bZ$ is $1$. Fix a $\bZ$-basis $(e_j)_{1\leq j\leq r}$ of $Q\cap L$, where $r=\dim_{\bQ}(Q)$. For $x\in P/Q$ and $n\in \bZ^r$, let $\sig(x, n)$ be the cone generated by all $N\in \fg$ such that the restriction of $N$ to $H'_\bR$ coincides with $N'$ and such that $N(e) =s(x) + a(x)^{-1}\sum_{j=1}^r t_je_j$ with $n_j \leq t_j \leq n_j+1$ for $1\leq j\leq r$.

Let $\Sig$ be the set of all faces of $\sig(x,n)$
for all $x\in P/Q$ and $n\in\bZ^r$.

\bigskip

\proclaim{Proposition 11.2} Let $\Gamma$ be as in \S$4$. 
Then
$\Sig$ is strongly
compatible with $\G$. 
\endproclaim

\noindent 
{\it Proof.} It is enough to 
show that $\Ad(\g)\sig\in \Sig$ for any $\g\in \G$ and $\sig\in \Sig$. It is sufficient to prove $\Ad(\g)\sig(x, n)\in \Sig$ for any $x\in P/Q$ and $n\in \bZ^r$. 
 Write 
$\g^{-1}e= e+h$ with $h\in H'$, and write
$\g s(x)+ N'(\g h)= s(y)+ q$ with $y\in P/Q$ and $q\in Q$. We have $a(y)=a(x)$ and 
 $a(x)q\in Q\cap L$. Write $a(x)q= \sum_{j=1}^r m_je_j$ with $m_j\in \bZ$, and let $m=(m_j)_j\in \bZ^r$. 
  We prove
$$\Ad(\g)\sig(x, n) = \sig(y, n+m).$$ In fact, let
$N\in \sig(x, n)$ and assume that the restriction of $N$ to $H'_{\bR}$ is $N'$ and 
$N(e)=s(x)+ a(x)^{-1}\sum_{j=1}^r t_je_j$  ($n_j\leq t_j \leq n_j+1$). 
Then, since $\g$ and $N'$ commute, we have 
$$\g N\g^{-1}(e)=\g N(e+h) = \g s(x)+ N'(\g h)+ a(x)^{-1}\sum_{j=1}^r t_je_j
$$ $$= s(y)+ a(y)^{-1}\sum_{j=1}^r (t_j + m_j)e_j.$$
  Hence $\g N\g^{-1}$ belongs to $\sig(y, n+m)$. \qed

\proclaim{Proposition 11.3}  $\Sig$ is a relatively complete fan.
\endproclaim

This is deduced from the 
following two facts:

\medskip

{\it Fact} 1.  Let $N :
H_\bR
\to H_\bR$ be a homomorphism such that $N(H_\bR)
\subset H'_\bR$ and such that the restriction of
$N$ to
$H'_\bR$ coincides with $N'$. Then the relative monodromy
filtration $M(N, W)$ exists if and only if
$N(e) \in \Im(N':H_\bR\to H_\bR)+W'_{-2,\bR}$.

\medskip

{\it Fact} 2. Let $N_j : H_\bR\to H_\bR$ ($j=1, 2$)  be 
homomorphisms such that $N_j(H_\bR)
\subset H'_\bR$ and such that the restrictions of
$N_j$ to
$H'_\bR$ coincide with $N'$. Then $N_1N_2=N_2N_1$ if and
only if
$N_1(e) \equiv N_2(e)\bmod \Ker(N' : H'_\bR\to
H'_\bR)$.

\head
\S12. Remarks on N\'eron models 
\endhead

  Here the base $S$ is any fs log analytic space unless otherwise stated. 
  Let $H'$ be a polarized log Hodge structure of weight $-1$ over $S$.  

\medskip

\noindent
{\bf {12.1.}}
  (The case where $F^1=0$ of this subparagraph is in [KKN08c] \S5.) 

  From the viewpoint of the theory of log intermediate Jacobian, it is fundamental 
to consider the exact sequences 
$0 \to H' \to (\cO^{\log}\otimes H')/F^0 \to
H'\bs(\cO^{\log} \otimes H')/F^0 \to 0$ 
of abelian sheaves on $(\fs/S)^{\log}$ (see \cite{KKN08a} for the 
definition of $(\fs/S)^{\log}$) and the induced 
$$0 \to \tau_*H' \to (H'_{\cO}/F^0)^{\text{hor}} \to 
(\tau_*(H'\bs(\cO^{\log} \otimes H')/F^0))^{\text{hor}} 
\overset \partial \to \to R^1\tau_*H'\tag{$\ast$}$$ 
of abelian sheaves on $(\fs/S)$, 
where \lq\lq hor'' means the horizontal parts, 
i.e., the parts consisting of sections corresponding to 
pre-log mixed Hodge structures that satisfy the small 
Griffiths transversality.

  There are several important subgroups of 
$(\tau_*(H'\bs(\cO^{\log} \otimes H')/F^0))^{\text{hor}}$ 
which are, respectively,  the 
inverse images of some subgroups of the monodromy group 
$R^1\tau_*H'$ under the connecting homomorphism 
$\partial$ in the last exact sequence $(*)$.
  The sheaf $\cExt^1_{LMHS}(\bZ, H')$ is one of them, which is the 
inverse image of the \lq\lq admissible part'' of the monodromy group. 
  From the viewpoint of log geometry, it is this sheaf that should be called 
the \lq\lq log intermediate Jacobian'', and 
what have been called log intermediate Jacobians 
so far in this
article should be called 
\lq\lq models of the log intermediate Jacobian.''

  Note that this sheaf $\cExt^1_{LMHS}(\bZ, H')$ is a group object, and, in a sense, 
log smooth (even when $H'$ corresponds to some 
geometric object which 
degenerate in the usual sense), as so are log complex tori introduced 
in \cite{KKN08a}. 
  The authors expect that it would be possible to generalize the theory 
of log complex tori and the theory of 
their proper models developed in \cite{KKN08a} and \cite{KKN08c} 
to the log intermediate Jacobians.

\medskip

\noindent
{\bf {12.2.}}
  Let the situation be as in \S1. 
  Let $\Sig_1$ be 
the fan consisting of $\{0\}$ and the cones 
$(\bR_{\ge0})N$ for $N\in\fg_\bQ$ satisfying $N|_{H'}=N'$, 
$N(e)=N'(a)$ for some $a \in H'_{\bQ}$ such that $\g a-a \in H'$. 
  We define the N\'eron model as $J^r_{X/S,\Sig_1}$.
  This is a log manifold whose log structure is the inverse image of that 
of the base, and 
\lq\lq represents'' (in some suitable senses
\footnote{Precisely, we can say as follows. 
In general, 
$J^r_{X/S,\Sig}$ is defined as an object of $\cB(\log)$ (\cite{KU09} 3.2.4)  
(even if $S$ is not $\Delta$ but any fs log 
analytic space). 
  If we do 12.1 over $\cB(\log)/S$ instead of $(\fs/S)$, 
we can say $J^r_{X/S,\Sig_1}$ represents the indicated subgroup 
in the usual sense. 
  Alternatively, in general, if the base $S$ is log smooth (\cite{KU09} 
2.1.11), 
$J^r_{X/S,\Sig}$ is a log manifold and so is regarded as an object of 
the category 
$(\fs/S)^{\sim}$ of sheaves on $(\fs/S)$. 
  Then, $J^r_{X/S,\Sig_1}$ is naturally isomorphic to 
the indicated subgroup as an object of $(\fs/S)^{\sim}$.
}) 
the subgroup of 
$(\tau_*(H'\bs(\cO^{\log} \otimes H')/F^0))^{\text{hor}}$ 
which is 
the inverse image of $\iota^{-1}((R^1 \tau_{S*} H')_{\text{tor}}) 
\subset R^1\tau_*H'$ by 
$\partial$ in $(*)$. 
  Here $\tau_S :S^{\log} \to S$ and $\iota:
(\fs/S) \to S$ are the natural morphisms. 
  Note that this $J^r_{X/S,\Sig_1}$ generalizes Clemens' model constructed 
in \S8 under the condition (8.1). 

  It is easy to see that for a sufficiently large $L$ (and for any $s$), 
our fan $\Sig$ in \S5 contains $\Sig_1$ as a subfan. 
  Hence our log intermediate Jacobian 
associated to $\Sig$ 
contains the N\'eron model as an open 
subspace. 
  In this sense, our construction gives a kind of 
compactification of the N\'eron model. 
  See \S8 for the special case of this fact under the condition (8.1). 

  By the proof of \cite{KU09} 4.3.1 (i), which works also in this 
mixed Hodge theoretic situation, 
any 
admissible normal function extends to the N\'eron model.

  The relationship does not seem to be known 
between this $J^r_{X/S,\Sig_1}$
and the N\'eron model constructed 
by Green-Griffiths-Kerr \cite{GGK.p}.

\enddocument